\documentclass[a4paper]{amsart}
\usepackage{a4wide,amssymb,amsmath,amsthm}

\vspace{2pt}

\theoremstyle{plain}
  \newtheorem{Th}{Theorem}
  \newtheorem*{Th*}{Theorem}
  \newtheorem{lm}{Lemma}
  
  \newtheorem{Pro}{Proposition}
\theoremstyle{definition}
  
  \newtheorem{defi}{Definition}
\theoremstyle{remark}
  
\theoremstyle{remark}
  
\newcommand{\R}{\mathbb{R}}
\newcommand{\g}{\widetilde{g}}
\newcommand{\MMg}{(\widetilde{M}^{n + m},\widetilde{g})}
\newcommand{\Mg}{(M^n,g)}
\newcommand{\MM}{\widetilde{M}^{n + m}}
\newcommand{\M}{M^n}
\newcommand{\MNg}{(\widetilde{M}^{n +1},\widetilde{g})}
\newcommand{\MN}{\widetilde{M}^{n + 1}}

\newcommand{\SH}{\mathbb{S}^k\times\mathbb{H}^{n - k + 1}}
\newcommand{\LV}{\widetilde{\nabla}}
\newcommand{\LVM}{\nabla}

\begin{document}

\title[hypersurfaces of $\SH$]{Hypersurfaces $\M$ in $\SH$}

\author[D. Kowalczyk]{Daniel Kowalczyk}
\address{Katholieke Universiteit Leuven\\ Departement
Wiskunde\\ Celestijnenlaan 200 B\\ B-3001 Leuven\\ Belgium}
\email[D. Kowalczyk]{daniel.kowalczyk@wis.kuleuven.be}

\begin{abstract}
Let $\psi:\M \rightarrow \SH$ be an isometric immersion of codimension $1$, then there exist symmetric $(1,1)$-tensors $S$ and $f$, a tangent vector field $U$ and a smooth function $\lambda$ on $\M$ that satisfy the compatibility equations of $\SH$. In this paper, we will deal with the converse problem: "Given a Riemannian manifold $\M$ with symmetric $(1,1)$-tensors $S$ and $f$, tangent vector field $U$ and smooth function $\lambda$ satisfying the conditions mentioned above, can $\M$ then be isometrically immersed in $\SH$ in such a way that $(g,S,f,U,\lambda)$ is realized as the induced structure?".
\end{abstract}


\subjclass[2000]{}

\maketitle
\section{Introduction}
It is well known that the Gauss and Codazzi equations are necessary conditions for a Riemannian manifold $\M$ to be locally isometrically immersed as hypersurface into an arbitrary Riemannian manifold $\MN$. In the case that $\MN(c)$ is a space form, the Gauss and Codazzi equations can be written in terms of the metric of $\M$ and of the shape operator $S$, which are both intrinsically known(as soon $S$ is known). Moreover in this case the Gauss and Codazzi equations are sufficient for an $n$-dimensional Riemannian manifold to be isometrically immersed into the space $\MN(c)$ as hypersurface with given second fundamental form $S$. In \cite{Da}, B. Daniel gave a necessary and sufficient condition for an n-dimensional Riemannian manifold to be isometrically immersed in the Riemannian product of a sphere and the real line or the Riemannian product of a hyperbolic space and the real line in terms of its first and second fundamental forms and of the projection of the vertical vector field $\partial_t$ on its tangent space. He rewrote the Gauss and Codazzi equations for hypersurfaces immersed in $\mathbb{S}^n\times \mathbb{E}^1$ or $\mathbb{H}^n \times \mathbb{E}^1$ in terms of the metric of $\M$, the shape operator $S$, the projection $T$ of the vertical vector field $\partial_t$ onto the tangent space of $\M$ and the normal component $\nu$ of $\partial_t$. The Gauss and Codazzi equations together with extra conditions on $T$ and $\nu$, which follow from the fact that $\partial_t$ is parallel in $\mathbb{S}^n\times \mathbb{E}^1$ or $\mathbb{H}^n\times \mathbb{E}^1$, are necessary and sufficient conditions for a Riemannian manifold to be isometrically immersed into $\mathbb{S}^n\times \mathbb{E}^1$ or into $\mathbb{H}^n\times \mathbb{E}^1$.
In this paper we will extend this result to the Riemannian product of the $k$-dimensional sphere and the $(n - k + 1)$-dimensional hyperbolic space of opposite sectional curvature.

On an arbitrary Riemannian product $M_1 \times M_2$ of two Riemannian manifolds there exist a natural symmetric $(1,1)$-tensor $F$ such that $F^2 = I\,(F \neq \pm I)$ and $\widetilde{\nabla}F = 0$, where $\widetilde{\nabla}$ is the Levi-Civita connection of the Riemannian product. $F$ is called the product structure of $M_1 \times M_2$. Let $\M$ be a hypersurface of the Riemannian product $M_1 \times M_2$ with product structure $F$. We can put
\[FX = fX + u(X)\xi,\]
\[F\xi = U + \lambda\xi,\]
where $X$ is a tangent vector field on $\M$, $\xi$ a unit normal of $\M$ in $M_1 \times M_2$, $f$ is a $(1,1)$-tensor on $M$, $U$ is a tangent vector field on $\M$, $u$ is a $1$-form on $\M$ and $\lambda$ is a smooth function on $\M$. Moreover there are some conditions on $f, U, u$ and $\lambda$ which follow from the fact that $F^2 = I$, $F$ is symmetric and $\widetilde{\nabla}F = 0$. In the case of $\mathbb{S}^k \times \mathbb{H}^{n - k + 1}$ one can rewrite the Gauss equation in terms of $S$ and $f$ and the Codazzi equation in terms of $u$, as follows:
\begin{gather*}
R(X,Y)Z = (SX \wedge SY)Z + \frac{1}{2}\left(f(X \wedge Y)Z + (X \wedge Y)fZ\right),\\
(\nabla_X S)(Y) - (\nabla_Y S)(X) = \frac{1}{2}\left(u(X)Y - u(Y)X\right).
\end{gather*}
Our main aim is to proof the following:
\begin{Th*}
Let $\Mg$ be a simply connected Riemannian manifold of dimension $n$. Let $S$ and $f$ be symmetric fields of operators on $\M$, $U$ a vector field on $\M$ and $\lambda$ a smooth function on $\M$ such that $f^2X = X - g(U,X)U, fU = -\lambda U$ and $\|U\|^2 + \lambda^2 = 1$. Assume that the structure $(g,S,f,U,\lambda)$ satisfies the Gauss and Codazzi equation for hypersurfaces in $\SH$ and the following equations:
\[(\nabla_X f)(Y) = g(U,Y)SX + g(SX,Y)U,\qquad \nabla_X U = \lambda SX - fSX,\qquad X[\lambda] = -2g(SX,U).\]
Then there exist an isometric immersion $\psi: \M \rightarrow \SH$ such that the shape operator with respect to the normal $\xi$ associated to $f$ is $S$ and such that
\begin{gather}
F(\xi) = \psi_{*}(U) + \lambda\xi,\\
F(\psi_{*}X) = \psi_{*}(f(X)) + g(U,X)\xi
\end{gather}
for every tangent vector field $X$ on $\M$ and whereby $F$ denotes the product structure of $\SH$. Moreover the immersion is unique up to the global isometries of $\SH$.
\end{Th*}
First we will recall some basic facts about the theory of submanifolds (\cite{C}) and some facts about Riemannian products. Then we will apply these facts to the case of hypersurfaces in the Riemannian product of the $k$-dimensional sphere and the $(n - k + 1)$-dimensional hyperbolic space. In the last section we will give the proof of the above stated theorem.
\section{Preliminaries}
\subsection{Submanifold Theory}
Let $\psi: \Mg \rightarrow \MMg$ be an isometric immersions of codimension $m$. Denote by $\nabla$ and $\widetilde{\nabla}$ the Levi-Civita connections of respectively $\Mg$ and $\MMg$. The tangent vector fields on $\M$ will be denoted by $X,Y,Z,\dots$ and normal vector fields along $\M$ in $\MM$ will be denoted by $\xi, \xi_1, \xi_2,\dots$. The formulas of Gauss and Weingarten which concern the decomposition of $\widetilde{\nabla}_X Y$ and $\widetilde{\nabla}_X\xi$, respectively, into their tangent and normal components are given by
\begin{gather*}
\LV_X Y = \LVM_X Y + h(X,Y),\\
\LV_X\xi = -S_{\xi}X + \nabla^{\perp}_X\xi,
\end{gather*}
whereby $h$ denotes the second fundamental form, $S_{\xi}$ denotes the shape operator of $\M$ in $\MM$ associated to the normal $\xi$ and $\nabla^{\perp}$ denotes the normal connection of $\M$ in $\MM$. The second fundamental form $h$ and the shape operator $S_{\xi}$ are related by
\[\widetilde{g}(h(X,Y),\xi) = g(S_{\xi}X,Y).\]
The shape operator $S_{\xi}$ at a point $p$ of $\M$ is a symmetric linear endomorphism on the tangent space $T_p\M$ for every normal vector $\xi$, since the second fundamental form $h$ is a normal-bundle-valued symmetric $(1,2)-$tensor field on $\M$.

Let $R,\widetilde{R}$ and $R^{\perp}$ denote the Riemannian curvature tensors of $\nabla,\LV$ and $\nabla^{\perp}$, respectively, taken with sign convention $R(X,Y) = [\LVM_X,\LVM_Y] - \LVM_{[X,Y]}$. One can easily deduce, using the formulae of Gauss and Weingarten that
\[\widetilde{R}(X,Y)Z = R(X,Y)Z + A_{h(X,Z)}Y - A_{h(Y,Z)}X + (\bar{\nabla}_Xh)(Y,Z) - (\bar{\nabla}_Yh)(X,Z),\]
where $\bar{\nabla}$ is the van der Waerden-Bortolotti connection of $\M$ in $\MM$ defined by
\[(\bar{\nabla}_X h)(Y,Z) = \nabla^{\perp}_X h(Y,Z) - h(\nabla_X Y,Z) - h(Y, \nabla_X Z).\]
The decomposition of $\widetilde{R}(X,Y)Z$ into a tangential component and normal component yields the equations of Gauss and Codazzi:
\begin{gather*}
\mathrm{tan}(\widetilde{R}(X,Y)Z) = R(X,Y)Z + A_{h(X,Z)}Y - A_{h(Y,Z)}X,\\
\mathrm{nor}(\widetilde{R}(X,Y)Z = (\bar{\nabla}_Xh)(Y,Z) - (\bar{\nabla}_Yh)(X,Z).
\end{gather*}
It is also possible to deduce the following equation of Ricci, again using the formulae of Gauss and Weingarten:
\[\mathrm{nor}(\widetilde{R}(X,Y)\xi) = R^{\perp}(X,Y)\xi + h(X,A_\xi Y) - h(A_\xi X,Y).\]

In the following we will mostly consider isometric immersions $\psi: \Mg \rightarrow \MNg$ of codimension 1, i.e. hypersurfaces. Let $\xi$ be a local unit normal vector field of $\M$ in $\MN$; the formulae of Gauss and Weingarten can be written as
\begin{gather*}
\LV_X Y = \LVM_X Y + h(X,Y)\xi,\\
\LV_X\xi = -SX,
\end{gather*}
whereby $h$ denotes the scalar valued second fundamental form of $\M$ in $\MN$ and $S$ denotes the shape operator of $\M$ in $\MN$, such that $g(SX,Y) = h(X,Y)$. The equations of Gauss and Codazzi then reduce to
\begin{gather*}
\mathrm{tan}(\widetilde{R}(X,Y)Z) = R(X,Y)Z - (SX \wedge SY)Z,\\
\mathrm{tan}(\widetilde{R}(X,Y)\xi) = (\nabla_X S)(Y) - (\nabla_Y S)(X),
\end{gather*}
whereby $\wedge$ associates to two tangent vectors $v,w \in T_p\M$ the endomorphism $v \wedge w$ of $T_p\M$ defined by
\[(v \wedge w)z = g(w,z)v - g(v,z)w\]
and whereby $\nabla S$ denotes the covariant derivative of the shape operator, i.e. whereby
\[(\nabla_X S)(Y) = \nabla_X(SY) - S(\nabla_X Y).\]
\subsection{Riemannian products}\label{prod}
Let $M_1$ and $M_2$ be two differentiable manifolds of dimensions $p \geq 2$ and $q \geq 2$, respectively. Consider the product manifold $M_1\times M_2$ with natural projections $\pi_1: M_1 \times M_2 \rightarrow M_1$ and $\pi_2: M_1 \times M_2 \rightarrow M_2$. The tangent space $T_{(p_1,p_2)}(M_1 \times M_2)$ of $M_1 \times M_2$ at $(p_1,p_2)$ is isomorphic to the direct sum $T_{p_1}M_1 \bigoplus T_{p_2}M_2$ for every point $(p_1,p_2)$ of $M_1 \times M_2$, i.e. we can identify the vector space $T_{(p_1,p_2)}(M_1 \times M_2)$ with the vector space $T_{p_1}M_1 \bigoplus T_{p_2}M_2$. At every point $(p_1,p_2)$ one can define the linear endomorphism $F_{(p_1,p_2)}:T_{(p_1,p_2)}(M_1 \times M_2) \rightarrow T_{(p_1,p_2)}(M_1 \times M_2): (v_1,v_2) \mapsto (v_1,-v_2)$. Since $F_{(p_1,p_2)}$ is thus defined at every point of $(p_1,p_2)$, we can construct a field of endomorphisms of the tangent spaces to $M_1 \times M_2$ by
\[F(X_1,X_2) = (X_1,-X_2),\]
whereby $X_1 \in TM_1$ and $X_2 \in TM_2$. Notice that $F^2 = I$, whereby $I$ denotes the identity transformation on $T(M_1 \times M_2)$. It is easy to see that a tangent vector $v$ to $M_1 \times M_2$ at $(p_1,p_2)$ lies in $T_{p_1}M_1$ if and only if $Fv = v$ and that a tangent vector $w$ to  $M_1 \times M_2$ at $(p_1,p_2)$ lies in $T_{p_2}M_2$ if and only if $Fw = -w$.

If $(M_1,g_1)$ and $(M_2,g_2)$ are Riemannian manifolds, consider the product Riemannian metric $\widetilde{g}$ on $M_1 \times M_2$:
\[\widetilde{g}((X_1,X_2),(Y_1,Y_2)) = g_1(X_1,Y_1) + g_2(X_2,Y_2),\]
for all $X_1,Y_1 \in TM_1$ and for all $X_2,Y_2 \in TM_2$. A product manifold of Riemannian manifolds equipped with the product Riemannian metric is called a Riemannian product manifold. Notice that the subspaces $T_{p_1}M_1$ and $T_{p_2}M_2$ of $T_{(p_1,p_2)}(M_1\times M_2)$ are orthogonal and that $F$ is a symmetric $(1,1)$-tensor on $M_1 \times M_2$ with respect to $\widetilde{g}$. Denote by $\LV$ the Levi-Civita connection of $(M_1 \times M_2,\widetilde{g})$. It can be shown that
\begin{equation}\label{par}
(\LV_X F)(Y) = 0,
\end{equation}
i.e. that
\[
\LV_X (FY) = F(\LV_X Y).
\]
Consider an isometric immersion $\psi: (\M,g) \rightarrow (M_1 \times M_2,\widetilde{g})$ of codimension $1$ with unit normal $\xi$, i.e. a hypersurface $\M$ in $M_1 \times M_2$. Then, we can put
\begin{gather}
FX = fX + u(X)\xi,\label{dec1}\\
F\xi = U + \lambda\xi,\label{dec2}
\end{gather}
whereby $f$ is a $(1,1)$-tensor on $\M$, $u$ is a $1$-form on $\M$, $U$ is a tangent vector field on $\M$ and $\lambda$ is a smooth function on $\M$.
\begin{lm}\label{sym1}
Suppose that $\M$ is a hypersurface of a Riemannian product $M_1 \times M_2$ with unit normal $\xi$ and for which $f, u, U$ and $\lambda$ are defined as above. Then
\begin{enumerate}
  \item f is a symmetric $(1,1)$-tensor such that $f^2X = X - u(X)U$ and that $fU = -\lambda U$,
  \item $u(X) = g(U,X)$, and
  \item $g(U,U) + \lambda^2 = 1.$
\end{enumerate}
\end{lm}

\begin{proof}
Using the fact that $F$ is a symmetric $(1,1)$-tensor field on $M_1 \times M_2$ with respect to the product Riemannian metric $\widetilde{g}$, we obtain that $f$ is a symmetric $(1,1)$-tensor field on $\M$ and that $u(X) = g(U,X)$ for every tangent vector field of $\M$:
\begin{gather*}
g(fX,Y) = \widetilde{g}(FX,Y) = \widetilde{g}(X,FY) = g(X,fY),\\
u(X) = \widetilde{g}(FX,\xi) = \widetilde{g}(X,F\xi) = g(U,X).
\end{gather*}
Using the equation $F^2 = I$ and comparing the tangent and the normal components, we obtain the other relations:
\begin{equation*}
\begin{split}
X & = F^2X\\
  & = F(fX + g(U,X)\xi)\\
  & = f^2X + g(U,fX)\xi + g(U,X)(U + \lambda\xi)
\end{split}
\end{equation*}
and
\begin{equation*}
\begin{split}
\xi & = F^2\xi\\
    & = F(U + \lambda\xi)\\
    & = fU + g(U,U)\xi + \lambda(U + \lambda\xi).
\end{split}
\end{equation*}
\end{proof}
Notice that the third equation can also be obtained from $fU = -\lambda U$, when $U \neq 0 $ in every point $p$ of $\M$.
\begin{lm}\label{parF1}
Let $\psi: \M \rightarrow M_1 \times M_2$ be an isometric immersion of codimension $1$ with local unit normal $\xi$ in a Riemannian product and for which $f, U$ and $\lambda$ are defined as above. Then
\begin{enumerate}
  \item $(\LVM_X f)(Y) = u(Y)SX + g(SX,Y)U$,
  \item $\LVM_X U = \lambda SX - fSX$
  \item $X[\lambda] = -2u(SX)$.
\end{enumerate}
\end{lm}
\begin{proof}
We will use equation $(\ref{par})$ to determine the covariant derivatives of $f, U$ and $\lambda$. We will therefore find a expression of $(\widetilde{\nabla}_X F)(Y)$ in terms of $f, U$ and $\lambda$ using the formulas of Gauss and Weingarten:
\begin{equation*}
\begin{split}
0 & = \widetilde{\nabla}_X FY - F(\widetilde{\nabla}_X Y)\\
  & = \widetilde{\nabla}_X(fY + g(U,X)\xi) - F(\nabla_X Y + g(SX,Y)\xi), and\\
  & = (\nabla_X f)Y - g(U,Y)SX - g(SX,Y)U + (g(\nabla_X U,Y) + g(fSX,Y) - \lambda g(SX,Y))\xi\\
\end{split}
\end{equation*}
and
\begin{equation*}
\begin{split}
0 & = \widetilde{\nabla}_X F\xi - F(\widetilde{\nabla}_X \xi)\\
  & = \widetilde{\nabla}_X (U + \lambda\xi) + F(SX)\\
  & = \nabla_X U - \lambda SX  + fSX + (X[\lambda] + 2g(U,SX))\xi.
\end{split}
\end{equation*}
By comparing the tangent and normal components, we obtain the equations of Lemma $2$. We notice that the third equation can also be obtained form the third equation of Lemma 1 and the second equation of Lemma 2 when $\|U\| \neq 1$ for every point $p$ of $\M$
\end{proof}

\section{Isometric immersions into $\SH$}
Denote by $\mathbb{L}^{n + 3} = (\R^{n + 3}, \langle.,.\rangle = dx_1^2 + \dots + dx_{n + 2}^2 - dx_{n + 3}^2)$ the Minkowski space of dimension $n + 3$ and define $\SH$ as the following submanifold of $\mathbb{L}^{n + 3}$, equipped with the induced metric:
\[\SH = \{(x_1,\dots,x_{n + 2},x_{n + 3})| x_1^2 + \dots + x_{k + 1}^2 = 1,\,x_{k + 2}^2 + \dots + x_{n + 2} - x_{n + 3}^2 = -1,\,x_{n + 3} > 0\}.\]
Then $\SH$ is the Riemannian product of the $k$-dimensional sphere and the $(n - k + 1)$-dimensional hyperbolic space. The Riemannian product $\SH$ is a submanifold of $\mathbb{L}^{n + 3}$ with orthonormal normal vectors $\xi_1 = (x_1,\dots,x_{k + 1},0,\dots,0)$ and $\xi_2 = (0,\dots,0,x_{k + 2},\dots,x_{n + 3})$. Using the formula of Gauss for isometric immersions into semi-Riemannian manifolds, we find that the Levi-Civita connection $\widetilde{\nabla}$ of $\SH$ is given by the following expression in terms of the Levi-Civita connection $D$ of $\mathbb{L}^{n + 3}$ and the product structure $F$ of $\SH$:
\begin{equation}\label{LVSH}
\widetilde{\nabla}_X Y = D_X Y + \frac{1}{2}\langle X + FX,Y\rangle\xi_1 - \frac{1}{2}\langle X - FX,Y\rangle\xi_2.
\end{equation}
Consequently, we obtain that the Riemann-Christoffel curvature tensor $\widetilde{R}$ is given by
\begin{equation}\label{curv2}
\widetilde{R}(X,Y)Z = \frac{1}{2}\left(F(X \wedge Y)Z + (X \wedge Y)FZ\right).
\end{equation}
Consider a hypersurface $\M$ of $\SH$ with unit normal $\xi$ and shape operator $S$. Let $\nabla$ and $R$ be the Levi-Civita connection and the Riemann-Christoffel curvature  of $\M$. Let $f$ be a endomorphism on $T\M$, $u$ a $1$-form, $U$ a vector field tangent to $\M$ and $\lambda$ a smooth function on $\M$ defined as in section $\ref{prod}$. Since the product structure $F$ of $\SH$ is symmetric and satisfies  $F^2= I$, we obtain the equations of Lemma $\ref{sym1}$ for $f, u, U$ and $\lambda$. Using these notations and equation $(\ref{curv2})$, we obtain that the equations of Gauss and Codazzi reduce to
\begin{gather}
R(X,Y)Z = (SX \wedge SY)Z + \frac{1}{2}\left(f(X \wedge Y)Z + (X \wedge Y)fZ\right),\\
(\nabla_X S)(Y) - (\nabla_Y S)(X) = \frac{1}{2}\left(u(X)Y - u(Y)X\right).
\end{gather}
We also have that the product structure $F$ of $\SH$ is parallel and hence we also obtain the equations of Lemma $\ref{parF1}$ for $f, U$ and $\lambda$. The equations of Lemmas $\ref{sym1}$ and $\ref{parF1}$ together with the  equations of Gauss and Codazzi for hypersurfaces of $\SH$ are called the compatibility equations of $\SH$. In the following we will suppose that $f \neq \pm Id$, because otherwise we would have that $U = 0$, $\lambda = \pm 1$ and
\begin{gather*}
R(X,Y)Z = (SX \wedge SY)Z  \pm(X \wedge Y)Z,\\
(\nabla_X S)(Y) - (\nabla_Y S)(X) = 0,
\end{gather*}
and hence $\M$ would be an isometrically immersed into $\mathbb{S}^{n + 1}$ or into $\mathbb{H}^{n + 1}$.

\begin{defi}
Let $S$ and $f$ be symmetric fields of operators on $\M$, $U$ a vector field on $\M$ and $\lambda$ a smooth function on $\M$. We say that the structure $(g,S,f,U,\lambda)$ of a Riemannian manifold $(\M,g)$ satisfies the compatibility equations of $\SH$ if the following equations are satisfied:
\begin{gather}
g(fX,Y) = g(X,fY) \qquad f^2X = X - g(U,X)U\qquad fU = -\lambda U \qquad g(U,U) + \lambda^2 = 1,\label{sym}\\
R(X,Y)Z = (SX \wedge SY)Z + \frac{1}{2}\left((fX \wedge fY)Z + (X \wedge Y)fZ\right),\label{Gauss}\\
(\nabla_X S)(Y) - (\nabla_Y S)(X) = \frac{1}{2}\left(u(X)Y - u(Y)X\right),\label{Codazzi}\\
(\LVM_X f)(Y) = u(Y)SX + g(SX,Y)U,\label{covf}\\
\LVM_X U = \lambda SX - fSX,\label{covU}\\
X[\lambda] = -2u(SX).\label{covlambda}
\end{gather}
\end{defi}

We will show now that an arbitrary Riemannian manifold $\Mg$ with structure $(g,S,f,U,\lambda)$ that satisfies the compatibility equations of $\SH$ can be isometrically immersed in $\SH$ and moreover that the immersion is unique up to isometries of $\SH$.

\begin{Th}
Let $\Mg$ be a simply connected Riemannian manifold of dimension $n$. Let $S$ and $f$ be symmetric fields of operators on $\M$, $U$ a vector field on $\M$ and $\lambda$ a smooth function on $\M$. Assume that $(g,S,f,U,\lambda)$ satisfies the compatibility equations for $\SH$. Then there exists an isometric immersion $\psi: \M \rightarrow \SH$ such that the shape operator with respect to the normal $\xi$ associated to $\psi$ is $S$ and such that
\begin{gather}
F(\xi) = \psi_{*}(U) + \lambda\xi,\\
F(\psi_{*}X) = \psi_{*}(f(X)) + g(U,X)\xi
\end{gather}
for every tangent vector field $X$ on $\M$ and whereby $F$ is the product structure of $\SH$. Moreover the immersion is unique up to a global isometries of $\SH$.
\end{Th}

\section{Proof of the theorem}

In the following we will give the proof of the theorem. We will use the techniques of \cite{D} and \cite{DNV} to proof the theorem. Let $TM^n$ be the tangent bundle of $\M$. Suppose $N = \M \times \R^3$ is the trivial bundle equipped with the Minkowski metric $g_{-1}$. Denote by $B = T\M \oplus_{W} N$ the orthogonal Whitney sum of $T\M$ and the trivial bundle $N$. The metric on $B$ will be denoted by $\widetilde{g}$. Let $\xi, \xi_1$ and $\xi_2$ be an orthonormal frame in $N$ such that $g_{-1}(\xi,\xi) = g_{-1}(\xi_1,\xi_1) = -g_{-1}(\xi_2,\xi_2) = 1$. We define a connection $D$ in $B$ by
\begin{gather}
D_X Y = \nabla_X Y + g(SX,Y)\xi -\frac{1}{2}g(X + fX,Y)\xi_1 + \frac{1}{2}g(X - fX,Y)\xi_2,\\
D_X \xi = -SX -\frac{1}{2}g(U,X)(\xi_1 + \xi_2),\\
D_X \xi_1 = \frac{1}{2}(X + f(X) + g(U,X)\xi),\\
D_X \xi_2 = \frac{1}{2}(X - f(X) - g(U,X)\xi).
\end{gather}
It is easy verified that the connection $D$ on $B$ is compatible with the metric $\widetilde{g}$ on $B$, i.e.
\[X\g(\eta_1,\eta_2) = \g(D_X\eta_1,\eta_2) + \g(\eta_1,D_X\eta_2)\]
for every tangent vector field $X$ on $\M$ and $\eta_1,\eta_2 \in B$. The curvature tensor of $B$ with connection $D$ will be denoted by $R^{D}$. The curvature tensor $R^D: T\M \times T\M \times B \rightarrow B$ is a trilinear map over the module  $C^{\infty}(\M)$ of smooth functions on $\M$ defined by
\[R^D(X,Y)\eta = D_XD_Y\eta - D_YD_X\eta - D_{[X,Y]}\eta.\]
We will show that the connection $D$ on $B$ is flat, i.e. $R^D = 0$.
\begin{lm}
$R^D = 0$.
\end{lm}
\begin{proof}
We will only calculate $R^D(X,Y)Z$ for arbitrary vector fields $X,Y,Z$ on $\M$. We will obtain that $R^D(X,Y)Z = 0$, because of the equations $(\ref{Gauss}),(\ref{Codazzi})$ and $(\ref{covf})$. Using the definition of the connection $D$, we have
\begin{equation*}
\begin{split}
D_XD_Y Z & = D_X(\nabla_Y Z + g(SY,Z)\xi - \frac{1}{2}(g(Y + fY,Z))\xi_1 + \frac{1}{2}(g(Y - fY,Z))\xi_2)\\
         & = \nabla_X\nabla_Y Z + g(SX,\nabla_Y Z)\xi\\
         & - \frac{1}{2}(g(X + fX,\nabla_Y Z))\xi_1  +  \frac{1}{2}(g(X - fX,\nabla_Y Z))\xi_2\\
         &   + (g(\nabla_X SY,Z)+ g(SY,\nabla_X Z))\xi + g(SY,Z)(-SX - \frac{1}{2}g(U,X)(\xi_1 + \xi_2))\\
         &   -\frac{1}{2}(g(\nabla_XY,Z)+ g(Y,\nabla_XZ) + g(\nabla_XfY,Z) + g(fY,\nabla_XZ))\xi_1\\
         &   -\frac{1}{2}(g(Y + fY,Z))(\frac{1}{2}(X + fX + g(U,X)\xi))\\
         &   + \frac{1}{2}(g(\nabla_XY,Z)+ g(Y,\nabla_XZ) - g(\nabla_XfY,Z) - g(fY,\nabla_XZ))\xi_2\\
         &   + \frac{1}{2}(g(Y - fY,Z))(\frac{1}{2}(X - fX - g(U,X)\xi))\\
\end{split}
\end{equation*}
and a similar equation when $X$ and $Y$ are interchanged. We also have
\begin{multline*}
D_{[X,Y]}Z = \nabla_{[X,Y]}Z + g(S[X,Y],Z)\xi\\
- \frac{1}{2}(g([X,Y] + f[X,Y],Z))\xi_1 + \frac{1}{2}(g([X,Y] - f[X,Y],Z))\xi_2.
\end{multline*}
Hence we obtain that $R^D(X,Y)Z$ is given by
\begin{equation*}
\begin{split}
&R(X,Y)Z - (SX \wedge SY)Z - \frac{1}{2}(f(X\wedge Y)Z + (X \wedge Y)fZ)\\
&+ (g(\nabla_X SY,Z) - g(\nabla_Y SX,Z) - g(S[X,Y],Z) - \frac{1}{2}g(U,X)g(Y,Z) + \frac{1}{2}g(U,Y)g(X,Z))\xi\\
&- \frac{1}{2}(g(U,X)g(SY,Z) + g(\nabla_X fY,Z) - g(U,Y)g(SX,Z) - g(\nabla_YfX,Z) - g(f[X,Y],Z))\xi_1\\
&+  \frac{1}{2}(g(SX,Z)g(U,Y) - g(SY,Z)g(U,X) - g(\nabla_XfY,Z) + g(\nabla_YfX,Z) + g(f[X,Y],Z))\xi_2.
\end{split}
\end{equation*}
From equation $(\ref{Gauss}), (\ref{Codazzi})$ and $(\ref{covf})$ we can conclude that $R^D(X,Y)Z$ vanishes for all tangent vector fields $X,Y$ and $Z$. The cases $R^D(X,Y)\xi, R^{D}(X,Y)\xi_{1}$ and $R^{D}(X,Y)\xi_2$  can be treated analogously using equations $(\ref{Codazzi}),(\ref{covf}),(\ref{covU})$ and $(\ref{covlambda})$. Since $R^D$ is a trilinear map, we obtain that $R^{D}(X,Y)\eta = 0$ for every $\eta \in B$.
\end{proof}
We define now a bundle map F on $B$ by
\begin{gather}
FX = fX + g(U,X)\xi,\\
F\xi = U + \lambda\xi,\\
F\xi_1 = \xi_1,\\
F\xi_2 = -\xi_2.
\end{gather}
In the following two lemmas we will show that $F^2 = I$, $F$ is symmetric with respect to $\widetilde{g}$ and the covariant derivative of the bundle map $F$ on $B$ is $0$.

\begin{lm}\label{symlm}
$F^2\eta = \eta$ and $\widetilde{g}(F\eta_1,\eta_2) = \widetilde{g}(\eta_1,F\eta_2)$.
\end{lm}
\begin{proof}
This follows immediately form the definition of the bundle map $F$ and $(\ref{sym})$.
\end{proof}

\begin{lm}\label{parF}
$(D_X F)(\eta) = 0$ for every $X \in T\M$ and every $\eta \in B$.
\end{lm}
\begin{proof}
This follows immediately form the definition of the bundle map $F$ and equations $(\ref{covf}),(\ref{covU})$ $,(\ref{covlambda})$.
\end{proof}

Let $B_1$ and $B_{-1}$ be subsets of $B$ defined respectively by
\[\{\eta \in B\,|\, F\eta = \eta\}\]
and
\[\{\eta \in B\,|\, F\eta = -\eta\}.\]
In the next proposition we will show that there exist orthonormal parallel sections $\widetilde{\eta}_1,\dots,\widetilde{\eta}_{k + 1},$ $\bar{\eta}_{1},\dots,\bar{\eta}_{n - k + 2}$ with $k \in \{1,\dots,n\}$ such that $F\widetilde{\eta}_{i} = \widetilde{\eta}_{i}$ and $F\bar{\eta}_{\alpha} = -\bar{\eta}_{\alpha}$ for $i \in \{1,\dots,k + 1\}$ and $\alpha \in \{1,\dots,n - k + 2\}$.

\begin{Pro}\label{sections}
Let $\Mg$ be a Riemannian manifold with structure $(g,S,f,U,\lambda)$ that satisfies the compatibility equations of $\SH$. Let $B$ be the vector bundle over $\M$ as defined above and $F$ the bundle map of $B$ as defined above, then there exist orthonormal parallel sections $\widetilde{\eta}_1,\dots,\widetilde{\eta}_{k + 1},\bar{\eta}_{1},\dots,\bar{\eta}_{n - k + 2}$ with $k \in \{1,\dots,n\}$ such that $\widetilde{\eta}_1,\dots,\widetilde{\eta}_{k + 1} \in B_1$ and that $\bar{\eta}_{1},\dots,\bar{\eta}_{n - k + 2} \in B_{-1}$.
\end{Pro}

\begin{proof}
Let $p$ be a point of $\M$ such that $f_p \neq \pm Id$. By definition of the bundle map $F$, we obtain that $Fv$ lies in the subspace $V$ of $B_p$ spanned by $T_p\M$ and $\xi_p$ for every $v \in T_p\M$. Analogously we have that $F\xi_p$ lies in $V$. Since $F$ is symmetric with respect to $\widetilde{g}$, $\widetilde{g}$ is positive definite when restricted to $V$, $F^2 = Id$, $F(V) \subset V$ and $f_p \neq Id$, there exist an orthonormal basis $\{\widetilde{v}_1,\dots,\widetilde{v}_{k},\bar{v}_1,\dots,\bar{v}_{n - k + 1}\}$ with $k \in \{1,\dots,n\}$ such that $F\widetilde{v}_i = \widetilde{v}_i$ and $F\bar{v}_{\alpha} = -\bar{v}_{\alpha}$. Since $F\xi_1 = \xi_1$ and $F\xi_2 = -\xi_2$, we know that there exist an orthonormal basis $\{\widetilde{w}_1,\dots,\widetilde{w}_{k + 1},\dots,\bar{w}_1,$ $\dots,\bar{w}_{n - k + 2}\}$ of $B_p$, with $\widetilde{g}(\bar{w}_{n - k + 2},\bar{w}_{n - k + 2}) = -1$, such that $F\widetilde{w}_{i} = \widetilde{w}_{i}$ and $F\bar{w}_{\alpha} = -\bar{w}_{\alpha}$. Since $R^D = 0$ and $\M$ is simply connected, there exist parallel sections $\widetilde{\eta}_1,\dots,\widetilde{\eta}_{k + 1},\bar{\eta}_1,\dots,\bar{\eta}_{n - k + 2}$ on $B$ such that $\widetilde{\eta}_i(p) = \widetilde{w}_{i}$ and $\bar{\eta}_{\alpha}(p) = \bar{w}_{\alpha}$. Moreover since the connection $D$ is compatible with the metric $\widetilde{g}$, we have that $\widetilde{\eta}_1,\dots,\widetilde{\eta}_{k + 1},\bar{\eta}_1,\dots,\bar{\eta}_{n - k + 2}$ are parallel orthonormal sections on $B$. From lemma $\ref{parF}$, we obtain that $F\widetilde{\eta}_i = \widetilde{\eta}_i$ and $F\bar{\eta}_{\alpha} = -\bar{\eta}_{\alpha}$. This completes the proof.
\end{proof}

We are ready to proof theorem 1.

\begin{proof}[Proof of theorem 1]
In Proposition $\ref{sections}$ we showed that there exist parallel orthonormal sections on $B$ $\widetilde{\eta}_1,\dots,\widetilde{\eta}_{k + 1},\bar{\eta}_{1},\dots,\bar{\eta}_{n - k + 2}$ with $k \in \{1,\dots,n\}$ such that $\widetilde{\eta}_1,\dots,\widetilde{\eta}_{k + 1} \in B_1$ and that $\bar{\eta}_{1},\dots,\bar{\eta}_{n - k + 2} \in B_{-1}$ with $k \in \{1,\dots,n\}$. In Lemma $\ref{sym}$ we showed that $F^2 = I$ and that $F$ is symmetric with respect to $\widetilde{g}$. Since $\widetilde{\eta}_1,\dots,\bar{\eta}_{n - k + 2}$ are parallel orthonormal sections, we know that $\xi_1 =  \sum_{i = 1}^{k + 1}\widetilde{g}(\xi_1,\widetilde{\eta}_i)\widetilde{\eta}_i + \sum_{\alpha = 1}^{n - k + 2}\epsilon_{\alpha}\widetilde{g}(\xi_1,\bar{\eta}_{\alpha})\bar{\eta}_{\alpha}$, where $\epsilon_1 = \dots = \epsilon_{n - k + 1} = 1$ and $\epsilon_{n - k + 2} = -1$. We obtain that $\widetilde{g}(\bar{\eta}_{\alpha},\xi_1) = 0$ for $\alpha \in \{1,\dots,n - k + 2\}$, because $F$ is symmetric with respect to $\widetilde{g}$, $F\xi_1 = \xi_1$ and $F\bar{\eta}_{\alpha} = -\bar{\eta}_{\alpha}$ and hence $\xi_1 = \sum_{i = 1}^{k + 1}\widetilde{g}(\xi_1,\widetilde{\eta}_i)\widetilde{\eta}_i$. Analogously, we obtain that $\xi_2 = \sum_{\alpha = 1}^{n - k + 2}\epsilon_{\alpha}\widetilde{g}(\xi_2,\bar{\eta}_{\alpha})\bar{\eta}_{\alpha}$. Using the components of $\xi_1$ and $\xi_2$ with respect to $\widetilde{\eta}_1,\dots,\widetilde{\eta}_{k + 1},\bar{\eta}_{1},\dots,\bar{\eta}_{n - k + 2}$ will construct an isometric immersion $\psi: \M \rightarrow \SH$. Define $\psi: \M \rightarrow \mathbb{L}^{n + 3}$ by
\[p \mapsto (\widetilde{g}_p(\widetilde{\eta}_1,\xi_1),\dots,\widetilde{g}_p(\widetilde{\eta}_{k + 1},\xi_1),\widetilde{g}_p(\bar{\eta}_1,\xi_2),\dots,\widetilde{g}_p(\bar{\eta}_{n - k + 2},\xi_2)).\]
It is easy to see that $\sum_{i = 1}^{k + 1} \widetilde{g}(\widetilde{\eta}_i,\xi_1)^2 = 1$ and $\sum_{\alpha = 1}^{n - k + 2} \epsilon_{\alpha}\widetilde{g}(\bar{\eta}_{\alpha},\xi_2)^2 = -1$ and hence $\psi(\M) \subset \SH$. Next we will show that $\psi$ is an immersion. Let $p$ be an arbitrary point of $\M$ and $v$ an arbitrary tangent vector at $p$  such that $\psi_{*}v = 0$. We obtain that
\[0 = v[\widetilde{g}(\widetilde{\eta}_i,\xi_1)] = \widetilde{g}(\widetilde{\eta}_i,D_v\xi_1) = \widetilde{g}(\widetilde{\eta}_i,\frac{1}{2}(v + fv + g(U,v)\xi))\]
and
\[0 = \widetilde{g}(\bar{\eta}_{\alpha},\frac{1}{2}(v - fv - g(U,v)\xi)).\]
We also have that
\[\begin{split}
F(\frac{1}{2}(v + fv + u(v)\xi) & = \frac{1}{2}(fv + f^2 + u(v)U + (g(U, v +fv +\lambda v))\xi)\\
                                & = \frac{1}{2}(v + fv + u(v)\xi)
\end{split}\]
and hence
\begin{equation}\label{immer1}\frac{1}{2}(v + fv +u(v)\xi)  = 0.\end{equation}
Analogously, we obtain that $F(\frac{1}{2}(v - fv - u(v))\xi) = -\frac{1}{2}(v - fv -u(v)\xi))$ and hence
\begin{equation}\label{immer2}
\frac{1}{2}(v - fv - u(v)\xi) = 0.
\end{equation}
Summing up the equations $(\ref{immer1})$ and $(\ref{immer2})$, we obtain that $v  = 0$ and hence $\psi$ is an immersion, because $p$ and $v$ were arbitrary. Next we will show that $\psi$ is isometric, i.e. $g(v,w) = \langle\psi_{*}v,\psi_{*}w\rangle$, where $\langle.,.\rangle$ is the Lorentzian inproduct on $\R^{n + 3}$.
\[\begin{split}
\langle\psi_{*}v,\psi_{*}w\rangle & = \frac{1}{4}(\sum_{i = 1}^{k  + 1}\widetilde{g}(v + fv + u(v)\xi,\widetilde{\eta}_{i})\widetilde{g}(w + fw + u(w)\xi,\widetilde{\eta}_{i}) \\
 & + \sum_{\alpha = 1}^{n - k + 2}\epsilon_{\alpha}\widetilde{g}(v - fv - u(v)\xi,\bar{\eta}_{\alpha})\widetilde{g}(w - fw - u(w)\xi,\bar{\eta}_{\alpha}))\\
 & = \frac{1}{4}(g(v + fv + u(v)\xi,w + fw + u(w)\xi) + g(v - fv - u(v)\xi,w - fw - u(w)\xi))\\
 & = \frac{1}{4}(2g(v,w) + 2g(f^2v,w) 2 +  2u(v)u(w))\\
 & = \frac{1}{4}(2g(v,w) + 2g(v,w) - 2u(v)u(w) + 2u(v)u(w))\\
 & = g(v,w).
\end{split}
\]
We can conclude that $\psi$ is an isometric immersion of $\M$ into $\SH$. We will show now that $N = (\widetilde{g}(\widetilde{\eta}_1,\xi),\dots,\widetilde{g}(\bar{\eta}_{n - l + 2},\xi))$ is  a unit normal of $\psi(\M)$ in $\SH$ and moreover the shape operator of $\M$ associated to the normal $N$ is given by $S$. Let $p$ be an arbitrary point of $\M$ and $v \in T_p\M$:
\[\begin{split}
\langle\psi_{*}v,N\rangle &= \sum_{i = 1}^{k + 1}\widetilde{g}(\frac{1}{2}(v + fv + u(v)\xi),\widetilde{\eta}_i)\widetilde{g}(\xi,\widetilde{\eta}_i)\\
 & + \sum_{\alpha = 1}^{n - k + 2}\epsilon_{\alpha}\widetilde{g}(\frac{1}{2}(v - fv - u(v)),\bar{\eta}_{\alpha})\widetilde{g}(\xi,\bar{\eta}_{\alpha})\\
 & = \widetilde{g}(\frac{1}{2}(v + fv + u(v)\xi),\xi) + \widetilde{g}(\frac{1}{2}(v - fv - u(v)),\xi)\\
 & = \frac{1}{2}(u(v) - u(v))\\
 & = 0
\end{split}\]
and hence $N$ is a unit normal of $\psi(\M)$. But we also have that $\langle N,N_1\rangle = \langle N,N_2\rangle = 0$, where $N_1 = (\widetilde{g}(\widetilde{\eta}_1,\xi_1),\dots,\widetilde{g}(\widetilde{\eta}_{k + 1},\xi_1),0,\dots,0)$ and $N_2 = (0,\dots,0,\widetilde{g}(\bar{\eta}_1,\xi_2),\dots,\widetilde{g}(\bar{\eta}_{n - k + 2},\xi_2))$. Therefore $N$ is a unit normal of $\psi(\M)$ in $\SH$. Next we show that $S$ is indeed the shape operator of $\psi(\M)$ in $\SH$ with respect to $N$:
\[\begin{split}
\langle \psi_{*}v,\widetilde{\nabla}_{\psi_{*}v}N\rangle &= \langle \psi_{*}v,D_{\psi_{*}v}N\rangle\\
 & = \widetilde{g}(\frac{1}{2}(v + fv + u(v)\xi),-Sw)\\
 & + \widetilde{g}(\frac{1}{2}(v - fv - u(v)\xi),-Sw)\\
 & = g(v,-Sw).
\end{split}\]
Suppose now that $F$ is the product structure of $\SH$. It is easy to deduce by direct calculations that $F(N) = \psi_{*}(U) + \lambda N$ and that
$F(\psi_{*}(X)) = \psi_{*}(fX) + u(X)N$. Finally we will prove that the isometric immersion is unique up to a isometry of $\SH$. Let $\psi_1,\psi_2: \M \rightarrow \SH$ be two isometric immersions of $\M$ in $\SH$ with unit normals $N_1$ and $N_2$, respectively, such that
\begin{gather}
\widetilde{\nabla}_{\psi_{*j}(X)}N_j = -\psi_{*j}(SX)\\
F(\psi_{*j}(X)) = \psi_{*j}(fX) + u(X)N_j \label{tan1}\\
F(N_j) = \psi_{*j}(U) + \lambda N_j\label{nor1},
\end{gather}
with $j = 1,2$. We will search now for an isometry $\phi$ of $\SH$ such that $\phi \circ \psi_1 = \psi_2$. Define a map $C: \M \rightarrow O^{1}(n + 3)$ by
\begin{gather*}
C_p(\psi_{*1}(X)) = \psi_{*2}(X),\\
C_p(N_1) = N_2,\\
C_p(\widetilde{N}_1) = \widetilde{N}_{2},\\
C_p(\bar{N}_1) = \bar{N}_2,
\end{gather*}
for all $X \in T_p\M$ and whereby $\widetilde{N}_j = (\psi^1_j,\dots,\psi^{k + 1}_{j},0,\dots,0)$ and $\bar{N}_j = (0,\dots,0,\psi^{k + 2}_j,\dots,\psi^{n + 3}_j)$, where $\psi^i_j$ are the components of $\psi_j$ with respect to the standard basis of $\mathbb{R}^{n + 3}$. We have identified here $T_{\psi_1(p)}\R^{n + 3}$ and $T_{\psi_{2}(p)}\R^{n + 3}$ with $\R^{n + 3}$. We will show that $C$ is an constant map, by showing that $(D_XC)(V) = D_XC(V) - C(D_XV) = 0 $ for every $X \in T_p\M$ and every vector field $V$  along $\psi_1$. To prove $(D_XC)(V) = 0$, it is sufficient to consider two cases: $V$ is a tangent vector field along $\psi$ or $V$ is a normal vector field of $\M$ in $\mathbb{L}^{n + 3}$. Assume first that $V  = \psi_{*1}(Y)$, where $Y$ is a vector field of $\M$. Then
\[\begin{split}
(D_XC)(\psi_{*1}(Y))& = D_X C(\psi_{*1}(Y)) - C(D_{X}\psi_{*1}(Y))\\
 & = D_{X} \psi_{*2}(Y) - C(D_{X}\psi_{*1}(Y))\\
 & = \psi_{*2}(\nabla_X Y) + g(SX,Y)N_2\\
 & -\frac{1}{2}(g(X,Y + fY))\widetilde{N}_2 + \frac{1}{2}(g(X,Y - fY))\bar{N}_2\\
 & - C(\psi_{*1}(\nabla_X Y) + g(SX,Y)N_1\\
 & -\frac{1}{2}(g(X,Y + fY))\widetilde{N}_1 + \frac{1}{2}(g(X,Y - fY))\bar{N}_1)\\
 & = 0.
\end{split}\]
Assume now that $V = N_1$, then
\[\begin{split}
(D_XC)(N_1)& = D_X C(N_1) - C(D_{X}N_1)\\
 & = -\psi_{*2}(SX) - \frac{1}{2}u(X)(\widetilde{N}_2 + \bar{N}_2) + C(\psi_{*1}(SX) + \frac{1}{2}u(X)(\widetilde{N}_1 + \bar{N}_1))\\
 & = 0.
\end{split}\]
One can also prove that $(D_XC)(\widetilde{N}_1) = (D_XC)(\bar{N}_1) = 0$, by using $D_X \widetilde{N}_{j} = \frac{1}{2}(\psi_{*j}(X + fX) + u(X)N_j)$ and $D_X \bar{N}_{j} = \frac{1}{2}(\psi_{*j}(X - fX) - u(X)N_j)$. Hence we obtain that $C$ is a constant map and can be identified with $\phi \in O^{1}(n + 3)$. Since $\phi(\widetilde{N}_1) = \widetilde{N}_2$ and $\phi(\bar{N}_1) = \bar{N}_2$, we obtain that $\phi(\psi_1) = \psi_2$.
\end{proof}

\end{document}